\documentclass[12pt]{article}

\input{xypic.tex} \xyoption{all}
\usepackage{amsthm,amssymb,latexsym,amsmath}

 \newcommand{\cpx}{\mathbb{C}}
 
\newcommand{\p}[1]{{\mathbb{P}^{#1}}} \newcommand{\op}[1]{{\cal O}_{\mathbb{P}^{#1}}}
  \newcommand{\cV}{{\cal V}}
\newcommand{\im}{{\rm Im}~}

\newtheorem{theorem}{Theorem}

\newtheorem{proposition}[theorem]{Proposition}

\newtheorem*{definition}{{\bf Definition}}

\newtheorem*{claim}{{\bf Claim}}

\begin{document}

\title{Admissible sheaves on $\p3$}
\author{Marcos Jardim \\ IMECC - UNICAMP \\
Departamento de Matem\'atica \\ Caixa Postal 6065 \\
13083-970 Campinas-SP, Brazil }

\maketitle

\begin{abstract}
Admissible locally-free sheaves on $\p3$, also known in the literature as
mathematical instanton bundles, arise in twistor theory,
and are in 1-1 correspondence with instantons on $\mathbb{R}^4$. In this
paper, we study admissible sheaves on $\p3$ from the algebraic geometric
point of view. We discuss examples and compare the admissibility condition
with semistability and splitting type. 
\end{abstract}

\tableofcontents

\newpage
\baselineskip18pt


\section*{Introduction} \label{intro}

An intense interest on the construction and classification of locally-free
sheaves on the 3-dimensional complex projective space started on the late 70's,
when twistor theory yielded a 1-1 correspondence between instantons
(i.e. anti-self-dual connections of finite $L^2$ norm) on $\mathbb{R}^4$ and
certain holomorphic vector bundles on $\p3$; this is the celebrated Penrose-Ward
correspondence \cite{Ma,WW}. This fact was later used by Atiyah, Drinfeld,
Hitchin and Manin to construct and classify all instantons \cite{ADHM}.

Since then, many authors have studied the so-called {\em mathematical}
(or {\em complex}) {\em instanton bundles}, defined in the literature
as rank 2 locally-free sheaves $E$ on $\p3$ with $c_1(E)=c_3(E)=0$ and
$c_2(E)=c>0$ satisfying $H^0(\p3,E(-1))=H^1(\p3,E(-2))=0$; see \cite{CTT}
for a recent brief survey of this topic. These correspond to $SL(2,\cpx)$
instantons on $\mathbb{R}^4$ of charge $c$. The correct generalization for
higher rank sheaves is given by Manin and Drinfeld (see \cite{Ma}), and
leads us to the key definition of this paper:

\begin{definition}
An admissible sheaf on $\p3$ is a coherent sheaf $E$ satisfying:
$$ H^p(\p3,E(k))=0 ~~ {\it for} ~~ p\leq1 ~,~ p+k\leq-1,
~~ {\it and} ~~ p\geq2 ~,~ p+k\geq0 $$
\end{definition}

Admissible locally-free sheaves of rank $r$ and vanishing first Chern class
are in 1-1 correspondence with $SL(r,\cpx)$ instantons on $\mathbb{R}^4$
\cite{Ma}. In this paper, we study mainly torsion-free admissible sheaves with
vanishing first Chern class, which can be regarded as a generalization of instantons.
Our focus is on the algebraic geometric properties of such objects, like semistability
and splitting type.

The paper is organized as follows. In Section \ref{s1} we remark that admissible sheaves
are in 1-1 correspondence with certain monads, exploring a few properties and some
examples in Section \ref{ex}. We then discuss how admissibility and semistability
compare with one another in Section \ref{s2}, and conclude with an analysis of the
splitting type of torsion-free admissible sheaves with vanishing first Chern class in
the last section.

\bigskip

\paragraph{Acknowledgment.}
Some of the results presented here were obtained in joint work with
Igor Frenkel \cite{FJ2}; we thank him for his continued support. We
also thank the organizers and participants of the XVIII Brazilian
Algebra Meeting.


\section{Monads} \label{s1}

Let $X$ be a smooth projective variety. A {\em monad} on $X$ is a
sequence $V_{\bullet}$ of the following form:
\begin{equation} \label{m1}
\cV_{\bullet} ~ : ~ 0 \to V_{-1} \stackrel{\alpha}{\longrightarrow}
V_{0} \stackrel{\beta}{\longrightarrow} V_{1} \to 0
\end{equation}
which is exact on the first and last terms. Here, $V_k$ are
locally free sheaves on $X$. The sheaf $E=\ker\beta/\im\alpha$ is
called the cohomology of the monad $\cV_{\bullet}$, also denoted
by $H^1(\cV_{\bullet})$.

In this paper, we will focus on the so-called {\em special monads}
on $\p3$, which are of the form:
$$ 0\to  V\otimes\op3(-1) \stackrel{\alpha}{\longrightarrow}
W\otimes\op3 \stackrel{\beta}{\longrightarrow} V'\otimes\op3(1) \to 0 ~, $$
where $\alpha$ is injective and $\beta$ is surjective. The existence of
such objects has been completely classified by Floystad in \cite{F}; let
$v=\dim V$, $w=\dim W$ and $v'=\dim V'$.

\begin{theorem} \label{exist}
There exists a special monad on $\p3$ as above if and only if at least
one of the following conditions hold:
\begin{itemize}
\item $w \geq 2v'+ 2$ and $w \geq v+v'$;
\item $w\geq v + v' + 3$.
\end{itemize} \end{theorem}

Monads appeared in a wide variety of contexts within algebraic geometry,
like the construction of locally free sheaves on complex projective spaces,
the study of curves in $\p3$ and surfaces in $\p4$. In this section, 
we will see how they are related to admissible sheaves on $\p3$. 

\begin{theorem} \label{T1}
Every admissible torsion-free sheaf $E$ on $\p3$ can be obtained
as the cohomology of a special monad
\begin{equation} \label{m2}
0\to  V\otimes\op3(-1) \stackrel{\alpha}{\longrightarrow}
W\otimes\op3 \stackrel{\beta}{\longrightarrow} V'\otimes\op3(1) \to 0 ~,
\end{equation}
where $V=H^1(\p3,E\otimes\Omega^2_{\p3}(1))$, 
$W=H^1(\p3,E\otimes\Omega^1_{\p3})$
and $V'=H^1(\p3,E(-1))$.
\end{theorem}
\begin{proof}
Manin proves the case $E$ being locally-free in \cite[p. 91]{Ma}, using
the Beilinson spectral sequence. However, the argument generalizes word by
word for $E$ being torsion-free; just note that the projection formula
$$ R^ip_{1*}\left(p_1^*\op3(k)\otimes p_2^* F\right)=
\op3(k)\otimes H^i(\p3,F) $$
holds for every torsion-free sheaf $F$, where $p_1$ and $p_2$ are the
natural projections of $\p3\times\p3$ onto the first and second factors.
\end{proof}

Clearly, the cohomology sheaf $E$ is always coherent, but more can
be said in particular situations. Note that $\alpha\in{\rm Hom}(V,W)\otimes\op1$
and $\beta\in{\rm Hom}(W,V')\otimes\op1$. Clearly, the surjectivity of
$\beta$ as a sheaf map implies that the localized map $\beta_x$ is
surjective for all $x\in\p3$, while the injectivity of $\alpha$
as a sheaf map implies that the localized map $\alpha_x$ is
injective only for generic $x\in\p3$.

\begin{theorem} \label{T2}
The cohomology $E$ of the monad
\begin{equation} \label{m4}
0 \to V\otimes\op2(-1) \stackrel{\alpha}{\longrightarrow}
W\otimes\op2 \stackrel{\beta}{\longrightarrow} V'\otimes\op2(1) \to 0
\end{equation}
is a coherent admissible sheaf with:
$$ {\rm rank}(E) = \dim W - \dim V - \dim V' ~~ , ~~ c_1(E) = \dim V' - \dim V $$
$$ ch_2(E) = \frac{1}{2}(\dim V + \dim V') ~~ {\rm and} ~~ ch_3(E) = \frac{1}{6}(\dim V - \dim V') ~ . $$
Moreover: 
\begin{itemize}
\item $E$ is torsion-free if and only if the localized maps $\alpha_x$
are injective away from a subset of codimension 2;
\item $E$ is reflexive if and only if the localized maps $\alpha_x$
are injective away from finitely many points;
\item $E$ is locally-free if and only if the localized maps
$\alpha_x$ are injective for all $x\in\p3$.
\end{itemize} \end{theorem}

\begin{proof}
The kernel sheaf ${\cal K}=\ker\beta$ is locally-free, and one has
the sequence:
\begin{equation} \label{ker}
0 \to V\otimes\op3(-1) \stackrel{\alpha}{\longrightarrow}
{\cal K} \to E \to 0
\end{equation}
so $E$ is clearly coherent. Notice also that:
$$ {\rm ch}(E) = \dim W - \dim V \cdot {\rm ch}(\op3(1))
   - \dim V' \cdot {\rm ch}(\op3(1)) $$
from which the calculation of the Chern classes of $E$ follows easily.

Taking the dual of the sequence (\ref{ker}), we obtain:
\begin{equation} \label{ker*}
0 \to E^* \to {\cal K}^* \stackrel{\alpha^*}{\longrightarrow} V^*\otimes\op3(1) 
\to {\rm Ext}^1(E,\op3) \to 0
\end{equation}
since ${\cal K}$ is locally-free. In particular,
${\rm Ext}^2(E,\op3)={\rm Ext}^3(E,\op3)=0$ and
$$ {\cal I} = {\rm supp}~{\rm Ext}^1(E,\op3) =
\{ x\in\p3 ~ | ~ \alpha_x ~ {\rm is~not~injective} ~ \} $$

So it is now enough to argue that $\cal C$ is torsion-free if and only if
$\dim I=1$ and that $\cal C$ is reflexive if and only if $\dim I=0$; the
third statement is clear.

Recall that the $m^{\rm th}$-singularity set of a coherent sheaf $\cal F$
is given by:
$$ S_m({\cal F}) = \{ X\in\p3 ~|~ dh({\cal F}_x) \geq 3-m \} $$
where $dh({\cal F}_x)$ stands for the homological dimension of
${\cal F}_x$ as an ${\cal O}_x$-module:
$$ dh({\cal F}_x) = d ~~~ \Longleftrightarrow ~~~
\left\{ \begin{array}{l}
{\rm Ext}^d_{{\cal O}_x}({\cal F}_x,{\cal O}_x) \neq 0 \\
{\rm Ext}^p_{{\cal O}_x}({\cal F}_x,{\cal O}_x) = 0 ~~ \forall p>d
\end{array} \right. $$

In the case at hand, we have that $dh({\cal F}_x) = 1$ if $X\in I$,
and $dh({\cal F}_x) = 0$ if $X\notin I$. Therefore
$S_0({\cal C})=S_1({\cal C})=\emptyset$, while $S_2({\cal C})=I$.
It follows that \cite[Proposition 1.20]{ST} :
\begin{itemize}
\item if $\dim I = 1$, then $\dim S_m({\cal C})\leq m-1$ for all $m<3$,
hence $\cal C$ is a locally 1$^{\rm st}$-syzygy sheaf;
\item  if $\dim I = 0$, then $\dim S_m({\cal C})\leq m-2$ for all $m<3$,
hence $\cal C$ is a locally 2$^{\rm nd}$-syzygy sheaf.
\end{itemize}
The desired statements follow from the observation that $\cal C$ is
torsion-free if and only if it is a locally 1$^{\rm st}$-syzygy sheaf, 
while $\cal C$ is reflexive if and only if it is a locally
2$^{\rm nd}$-syzygy sheaf \cite[p. 148-149]{OSS}.
\end{proof}

As part of the proof above, it is worth emphasizing that if $E$ is
admissible then ${\rm Ext}^2(E,\op3)={\rm Ext}^3(E,\op3)=0$.

It follows from Theorems \ref{T1} and \ref{T2} that there exists a 
(set theoretical) 1-1 correspondence between special monads and 
admissible sheaves. As shown by Manin and Drinfeld (see \cite{Ma}),
such correspondence is categorical.

\begin{theorem}
The functor that associates a special monad on $\p3$ to its cohomology sheaf
defines an equivalence between the categories of special monads and admissible
sheaves.
\end{theorem}

We complete this section with an important fact:

\begin{proposition}\label{dv}
If $E$ is an admissible sheaf, then $H^0(\p2,E^*(k))=0$ for all $k\leq-1$.
\end{proposition}
\begin{proof}
$E$ is the cohomology of the monad (\ref{m2}); setting
$V=H^1(\p2,E(-2))$, $W=H^1(\p2,E\otimes\Omega^1_{\p2}(-1))$ and
$V'=H^1(\p2,E(-1))$, one had the sequences 
$$ 0 \to {\cal K}(k) \to W\otimes\op2(k) \to V'\otimes\op2(k+1) \to 0 ~~ {\rm and} $$
$$ 0 \to V\otimes\op2(k-1) \to {\cal K}(k) \to E(k) \to 0 ~ . $$
where ${\cal K}=\ker\{W\otimes\op2\to V'\otimes\op2\}$ is a locally-free sheaf.
It follows from the first sequence that:
$$ H^0(\p2,{\cal K}(k))=0 ~~ \forall k\leq-1 ~,~
H^2(\p2,{\cal K}(k))=0 ~~ \forall k\geq-2 $$
$$ {\rm and} ~~ H^0(\p2,{\cal K}^*(k))=0 ~~ \forall ~ k\leq-1 ~~, ~~ {\rm by~Serre~duality}.$$
The proposition then follows easily from the dual of the second sequence. 
\end{proof}


\section{Examples of admissible sheaves}\label{ex}

Let us now study various examples of admissible sheaves on $\p3$.
Theorem \ref{exist} implies that there are admissible
coherent sheaves in rank 0 and 1, but there can be no admissible sheaves
with zero first Chern class in these ranks, apart from the trivial ones. 

Examples of admissible sheaves with vanishing first Chern class start
in rank 2. The basic one is an admissible torsion-free sheaf $E$ which is not
locally-free; it arises as the cohomology $E$ of the monad:
\begin{equation} \label{ex-tf}
\op3(-1) \stackrel{\alpha}{\rightarrow} \op3^{\oplus4}
\stackrel{\beta}{\rightarrow} \op3(1)
\end{equation}
$$ \alpha = \left(\begin{array}{c} x \\ y \\ 0 \\ 0 \end{array}\right) 
~~{\rm and}~~ \beta= (-y ~~ x ~~ z ~~ w) $$
It is easy to see that $\beta$ is surjective for all $[x:y:z:w]\in\p3$, 
while $\alpha$ is injective provided $x,y\neq0$. It then follows from
Theorem \ref{T2} that $E$ is torsion-free, but not locally-free. In particular,
the singularity set of $E$ (i.e. the support of $E^{**}/E$) consists of
the line $\{x=y=0\}\subset\p3$. Note also that $c_2(E)=1$ and $c_1(E)=c_3(E)=0$.

Reflexive sheaves on $\p3$ have been extensively studied in a series of 
papers by Hartshorne \cite{Ha}, among other authors. In particular, it was
show that a rank 2 reflexive sheaf ${\cal F}$ on $\p3$ is locally-free if and only
if $c_3({\cal F})=0$. Therefore, we conclude:

\begin{proposition} ({\bf Hartshorne \cite{Ha}})
There are no rank 2 admissible sheaves on $\p3$ which are reflexive but not 
locally-free.
\end{proposition}

The situation for higher rank is quite different, though, and it is easy to 
construct a rank 3 admissible sheaf which is reflexive but not locally-free.
Setting $w=5$ and $v=v'=1$, consider the monad:
\begin{equation} \label{ex-ref}
\op3(-1) \stackrel{\alpha}{\rightarrow} \op3^{\oplus5}
\stackrel{\beta}{\rightarrow} \op3(1)
\end{equation}
$$ \alpha = \left(\begin{array}{c} x \\ y \\ 0 \\ 0 \\ z \end{array}\right) 
~~{\rm and}~~
\beta= (-y ~~ x ~~ z ~~ w ~~ 0) $$
Again, it is easy to see that $\beta$ is surjective for all $[x:y:z:w]\in\p3$,
while $\alpha$ is injective provided $x,y,z\neq0$. It then follows from
Theorem \ref{T2} that $EE$ is reflexive, but not locally-free; its singularity
set is just the point $[0:0:0:1]\in\p3$. Note also that $c_2(E)=1$ and $c_1(E)=c_3(E)=0$.

Finally, we give an example of a rank 2 admissible locally-free sheaf.
Setting $w=4$ and $v=v'=1$, consider the monad:
\begin{equation} \label{ex-lf}
\op3(-1) \stackrel{\alpha}{\rightarrow} \op3^{\oplus4}
\stackrel{\beta}{\rightarrow} \op3(1)
\end{equation}
$$ \alpha = \left(\begin{array}{c} x \\ y \\ -w \\ z \end{array}\right) 
~~{\rm and}~~
\beta= (-y ~~ x ~~ z ~~ w) $$
It is easy to see that $\beta$ is surjective and $\alpha$ is injective
for all $[x:y:z:w]\in\p3$, so $E$ is indeed locally-free; note that
$c_2(E)=1$ and $c_1(E)=c_3(E)=0$.

With these simple examples in low rank, we can produce high rank 
admissible sheaves using the following:

\begin{proposition} \label{ext}
If $F'$ and $F''$ are coherent admissible sheaves, its extention $E$:
$$ 0 \to F' \to E \to F'' \to 0 $$
is also admissible.
\end{proposition}

The proof is an easy consequence of the associated long exact sequence
in cohomology, and it is left to the reader. As a consequence of Serre
duality, we have:

\begin{proposition} \label{tensor-dual}
If $E$ is a locally-free admissible sheaf, then $E^*$ is also admissible.
\end{proposition}


\section{Semistability of torsion-free admissible sheaves} \label{s2}

Recall that a torsion-free sheaf $E$ on $\p3$ is said to be
{\em semistable} if for every coherent subsheaf $0\neq F\hookrightarrow E$
we have
$$ \mu(F)=\frac{c_1(F)}{{\rm rk}(F)} \leq \frac{c_1(E)}{{\rm rk}(E)}=\mu(E) ~ . $$
Furthermore, if for every coherent subsheaf $0\neq F\hookrightarrow E$
with $0<{\rm rk}(F)<{\rm rk}(E)$ we have
$$ \frac{c_1(F)}{{\rm rk}(F)} < \frac{c_1(E)}{{\rm rk}(E)} ~ , $$
then $E$ is said to be {\em stable}. It is also important to
remember that:
\begin{itemize}
\item $E$ is (semi)stable if and only if $E^*$ is;
\item $E$ is (semi)stable if and only if $\mu(F)<\mu(E)$ ($\mu(F)\leq\mu(E)$)
for all coherent subsheaves $F\hookrightarrow E$ whose quotient $E/F$ is
torsion-free;
\item $E$ is (semi)stable if and only if $\mu(Q)>\mu(E)$ ($\mu(Q)\geq\mu(E)$)
for all torsion-free quotients $E\to Q\to 0$ with $0<{\rm rk}(Q)<{\rm rk}(E)$.
\end{itemize}
Furthermore, if $E$ is locally-free, it is enough to test the
locally-free subsheaves $F\hookrightarrow E$ with $0<{\rm rk}(F)<{\rm rk}(E)$
to conclude that $E$ is stable.


The goal of this section is to compare the semistability and admissibility conditions.
We provide to positive results for admissible sheaves of rank 2 and 3. 

\begin{theorem} \label{ss-adm}
Let $E$ be a semistable torsion-free sheaf with $c_1(E)=0$. $E$ is admissible if and only if
$H^1(\p3,E(-2))=H^2(\p3,E(-2))=0$. Furthermore, an admissible torsion-free
sheaf is stable if and only if $H^0(\p3,E)=0$.
\end{theorem}

In other words, if $E$ is a semistable torsion-free sheaf with $c_1(E)=0$ and
$H^1(\p3,E(-2))=H^2(\p3,E(-2))=0$, then $E$ is admissible. 

\begin{proof}
Semistability implies immediately that $H^0(\p3,E(k))=0$, for all $k\leq1$. If
$E$ is locally-free, then by Serre duality we have $H^3(\p3,E(k))=0$ for all
$k\geq-3$, since $E^*$ is also semistable. If $E$ is torsion-free, we can use
the semistability of $E^**$ and the sequence
$$ 0 \to E \to E^{**} \to Q \to 0 ~~,~~ Q=E^{**}/E $$
to conclude that $^3(\p3,E(k))=^3(\p3,E^{**}(k))=0$, since $Q$ is supported in
dimension less or equal to 1.

Now we assume that $H^1(\p3,E(-2))=H^2(\p3,E(-2))=0$, and let $\wp$ be a plane
in $\p3$. From he sequence:
\begin{equation} \label{wp}
0 \to E(k-1) \to E(k) \to E|_{\wp}(k) \to 0
\end{equation}
we conclude that $H^0(\wp,E|_{\wp}(-1))=H^2(\wp,E|_{\wp}(-2))=0$.

\begin{claim}
If $V$ is a torsion-free sheaf on $\p2$ with $H^0(\p2,V(-1))=H^2(\p2,V(-2))=0$,
then $H^0(\p2,V(k))=0$ for $k\leq-1$ and $H^2(\p2,V(k))=0$ for $k\geq-2$.
\end{claim}
\noindent {\it Proof of the claim:} For any line $\ell\subset\p2$, we have the
sequence
$$ 0 \to V(k-1) \to V(k) \to V|_{\ell}(k) \to 0 ~~, $$
so that
$$ 0 \to H^0(\p2,V(k-1)) \to H^0(\p2,V(k)) ~~ {\rm and} ~~
H^2(\p2,V(k-1)) \to H^2(\p2,V(k)) \to 0 ~~. $$
The claim follows easily by induction.

Returning to (\ref{wp}), we also have:
$$ H^0(\wp,E|_{\wp}(k)) \to H^1(\p3,E(k-1)) \to H^1(\p3,E(k)) $$
so, for $k\leq-1$, if $H^1(\p3,E(k))=0$, then also $H^1(\p3,E(k-1))=0$.
Thus by induction we conclude that $H^1(\p3,E(k))=0$ for all $k\leq-2$.

Similarly, we have:
$$ H^2(\p3,E(k-1)) \to H^2(\p3,E(k)) \to H^2(\wp,E|_{\wp}(k)) $$
and again by induction we conclude that $H^2(\p3,E(k))=0$ for all
$k\geq-2$, as desired.
\end{proof}

The converse statement seems to depend on the rank, as we will see in the two
results below.

\begin{theorem}\label{ss-rk2}
Every rank 2 admissible torsion-free sheaf $E$ with $c_1(E)=0$ is semistable.
Moreover, if $H^0(\p3,E^*)=0$, then $E$ is stable.
\end{theorem}
\begin{proof}
First, assume that $L$ is a rank 2 reflexive sheaf with $c_1(L)=0$ and $H^0(L(k))=0$
for all $k\leq-1$; we show that $L$ is semistable. Indeed, let
$F\hookrightarrow L$ be a torsion-free subsheaf of rank 1, with torsion-free
quotient $Q=L/F$. By Lemma 1.1.16 in
\cite[p. 158]{OSS}, it follows that $F$ is also reflexive; but every rank 1
reflexive sheaf is locally-free, thus $F=\op3(d)$. Any map $F\to L$ yields
a section in $H^0(\p3,L(-d))$, $c_1(F)=d\leq0$ and $E$ is semistable, being
stable if $H^0(\p3,E^*)=0$

Now if $E$ is a rank 2 admissible torsion-free sheaf with $c_1(E)=0$, then
$L=E^*$ is a rank 2 reflexive sheaf with $c_1(L)=0$ and $H^0(L(k))=0$
for all $k\leq-1$, by Proposition \ref{dv} and since the dual of any
coherent sheaf is always reflexive. Thus $E^*$ is semistable, so $E$ is
as well. Clearly, $E$ is stable if $H^0(\p3,E^*)=0$, as desired.
\end{proof}

A similar result for rank 3 sheaves requires a stronger hypothesis: reflexivity,
rather than torsion-freeness.

\begin{theorem}\label{ss-rk3}
Every rank 3 admissible reflexive sheaf $E$ with $c_1(E)=0$ is semistable.
Moreover, if $H^0(\p3,E)=H^0(\p3,E^*)=0$, then $E$ is stable.
\end{theorem}
\begin{proof}
In fact, one can show that every rank 3 reflexive sheaf with $c_1(E)=0$
and $H^0(E(k))=H^0(E^*(k))=0$ is semistable. The desired theorem follows
easily from this fact.

Indeed, let $F\hookrightarrow E$ be a torsion-free subsheaf, with torsion-free
quotient $Q=E/F$, so that $c_1(F)=-c_1(Q)$. As in the proof of Theorem
\ref{ss-rk2}, it follows that $F$ is reflexive. There are two possibilities:

\noindent {\bf (i)} ${\rm rank}~F=1$. 
In this case, $F$ is locally-free, so a map $F\to E$ yields a section in
$H^0(\p3,E(-d))$, where $d=c_1(F)$. Hence $c_1(F)\leq0$.

\noindent {\bf (ii)} ${\rm rank}~F=2$, so ${\rm rank}~Q=1$. 
Now $Q^*$ is a reflexive (hence locally-free) subsheaf of $E^*$, which
gives a section in $H^0(\p3,E^*(-d))$, where $d=c_1(Q^*)=c_1(F)$. Hence
$c_1(F)\leq0$.

It follows that $E$ is semistable, being stable if
$H^0(\p3,E)=H^0(\p3,E^*)=0$.
\end{proof}

Together with Theorem \ref{ss-adm}, we conclude that:
\begin{itemize}
\item a rank 2 torsion-free sheaf on $\p3$ with $c_1(E)=0$ is
admissible if and only if it is semistable and $H^1(\p3,E(-2))=H^2(\p3,E(-2))=0$;
\item a rank 3 reflexive sheaf on $\p3$ with $c_1(E)=0$ is
admissible if and only if it is semistable and $H^1(\p3,E(-2))=H^2(\p3,E(-2))=0$.
\end{itemize}


\section{Trivial splitting type} \label{s3}

Since every locally-free sheaf on a projective line splits as a sum of
line bundles, one can study sheaves on projective spaces by looking into
the behavior of restriction to a line \cite{OSS}.

\begin{definition}
A torsion-free sheaf $E$ on $\p3$ is said to be of trivial splitting type if
there is a line $\ell\subset\p3$ such that $E|_\ell$ is the
trivial locally-free sheaf, i.e. $E|_\ell\simeq{\cal O}_{\ell}^{\oplus{\rm rk}E}$.
\end{definition}

A sheaf of trivial splitting type necessarily has vanishing first Chern class.
Note that, by semicontinuity, if $E$ is of trivial splitting type then
$E|_\ell$ is trivial for a generic line in $\p3$. Torsion-free sheaves of
trivial splitting type where completely classified in \cite{FJ2}, and they
were shown to be closely related with a complex version of the celebrated
Atiyah-Drinfeld-Hitchin-Manin matrix equations.

Furthermore, every torsion-free
sheaf of trivial splitting type is semistable; indeed, assume that $E$
has rank $r$, and let $F\hookrightarrow E$ be a coherent subsheaf of
rank $s$, with torsion-free quotient $E/F$. Then on a generic line
$\ell\subset\p3$ we have:
$$ F_\ell = {\cal O}_{\ell}(a_1)\oplus\cdots\oplus{\cal O}_{\ell}(a_s)
\hookrightarrow E|_\ell\simeq{\cal O}_{\ell}^{\oplus r} ~~, $$
where $c_1(F)=a_1+\cdots+a_s$. It follows that $c_1(F)\leq0$, since we
must have $a_k\leq0$, $k=1,\dots,s$.

\begin{theorem} \label{tst-adm}
Let $E$ be a torsion-free sheaf of trivial splitting type. $E$ is admissible
if and only if $H^1(\p3,E(-2))=H^2(\p3,E(-2))=0$.
\end{theorem}

Of course, this is an easy consequence of Theorem \ref{ss-adm} and the
observation above, but here is a direct proof.

\begin{proof}
Let $E$ be an admissible torsion-free sheaf.
Without loss of generality, we can assume that $E|_{\ell_\infty}$ is trivial
for $\ell_\infty=\{z=w=0\}$.
Let $\wp$ be a plane containing $\ell_\infty$, e.g. $\wp=\{z=0\}$. Then
$E|_{\wp}$ is a torsion-free sheaf on $\wp$ which is trivial at $\ell_\infty$.
From the proof of Theorem \ref{ss-adm} we know that:
\begin{equation} \label{v1}
H^0(\wp,E|_{\wp}(k)) = 0 ~ \forall k\leq-1 ~~ , ~~
H^2(\wp,E|_{\wp}(k)) = 0 ~ \forall k\geq-2
\end{equation}
Now consider the sheaf sequence:
\begin{equation} \label{v2}
0 \to E(k-1) \stackrel{\cdot z}{\longrightarrow} E(k)
\longrightarrow E|_{\wp}(k) \to 0
\end{equation}
Using (\ref{v1}), we conclude that:
$$ H^3(\p3,E(k)) =  H^3(\p3,E(k-1)) ~ \forall k\geq-2 $$
But, by Serre's vanishing theorem, $H^3(\p3,E(N))=0$ for sufficiently 
large $N$, thus
$H^3(\p3,E(k)) = 0$ for all $k\geq-3$.

Similarly, we have:
$$ H^0(\p3,E(k-1)) =  H^0(\p3,E(k)) ~ \forall k\leq-1 $$
Since $E\hookrightarrow E^{**}$, we have via Serre duality:
$$ H^0(\p3,E(k))\hookrightarrow H^0(\p3,E^{**}(k)) = 
H^3(\p3,E^{***}(-k-4))^*~. $$
Thus, again by Serre's vanishing theorem, $H^0(\p3,E(-N))=0$ for
for sufficiently large $N$, so that $H^0(\p3,E(k)) = 0$ for all 
$k\leq-1$.

We also have that:
$$ 0 \to H^1(\p3,E(k-1)) \to  H^1(\p3,E(k)) ~ \forall k\leq-1 $$
hence $H^1(\p3,E(-2)) = 0$ implies that $H^1(\p3,E(k)) = 0$ for all 
$k\leq-2$. Furthermore,
$$ H^2(\p3,E(k-1)) \to  H^2(\p3,E(k)) \to 0 ~ \forall k\geq-2 $$
forces $H^2(\p3,E(k)) = 0$ for all $k\geq-2$ once $H^2(\p3,E(-2)) = 0$.
\end{proof}

As in the previous section, the converse statement seems to depend on
the rank. The generic splitting type of a semistable locally-free sheaf
with vanishing first Chern class is determined by Theorem 2.1.4 in
\cite[p. 205-206]{OSS}. In particular, it follows that every semistable
rank 2 locally-free sheaf is of trivial splitting type. Thus, from
Theorem \ref{ss-rk2}, we conclude:

\begin{theorem} \label{tst-rk2}
Every rank 2 admissible locally-free sheaf $E$ with $c_1(E)=0$ is of
trivial splitting type.
\end{theorem}

To explore a few easy consequences of the classical theory of locally-free
sheaves on complex projective spaces, let $\mathbb{G}$ denote the Grasmannian
of lines in $\p3$.

\begin{definition}
Let $E$ be a locally-free sheaf of trivial splitting type. The set
$$ {\cal J}_E =\{ \ell\in\mathbb{G} ~~|~~ E_{\ell} ~ {\rm is~not~trivial} \} $$
is called the set of jumping lines of $E$; it is always a closed subvariety
of $\mathbb{G}$. Moreover, $E$ is said to be uniform if ${\cal J}_E$ is empty,
i.e. if $E_{\ell}$ is independent of $\ell\in\mathbb{G}$.
\end{definition}

\begin{theorem}\label{uni}
Every rank 2 uniform, admissible locally-free sheaf $E$ with $c_1(E)=0$ is trivial.
\end{theorem}
\begin{proof}
By Theorem \ref{tst-rk2}, $E$ is of trivial splitting type. Since $E$ is uniform,
$E_\ell$ must be trivial for all $\ell\in\mathbb{G}$. It then follows from Theorem
3.2.1 in \cite[p. 51]{OSS} that $E$ is trivial.
\end{proof}

Our last result, regarding the set of jumping lines, follows from Theorem 2.2.3
in \cite[p. 228]{OSS}.

\begin{theorem}
If $E$ is an rank 2 admissible locally-free sheaf with $c_1(E)=0$, then its
set of jumping lines ${\cal J}_E$ is a subvariety of pure codimension 1 in $\mathbb{G}$.
\end{theorem}


\section{Open problems}

The results proved in this paper point to a number of quite interesting questions
and possible generalizations. First of all, we expect that if $E$ is a properly
torsion-free or properly reflexive admissible sheaf, then its dual $E^*$ is not
admissible, but we have not been able to construct any examples.

We would also like to see whether the results in Section \ref{s2} can be generalized
to higher rank. It seems too much to expect every admissible sheaf to be semistable;
but the correspondence between instantons and locally-free admissible sheaves
makes the statement "every admissible locally-free sheaf is semistable"
an attractive conjecture. On the other hand, is Theorem \ref{ss-rk3} optimal, i.e.
is there a rank 3 torsion-free admissible sheaf which is not semistable?

It would also be interesting to study the connection between admissibility and
Gieseker stability. Since every Gieseker semistable sheaf on a projective space
is also semistable, we conclude from Theorem \ref{ss-adm} that every Gieseker
semistable torsion-free sheaf is admissible; one would like to determine to
what extent the converse is also true.

Theorems \ref{tst-rk2} and \ref{uni} point to interesting properties of higher
rank admissible sheaves: is every admissible locally-free sheaf with vanishing
first Chern class of trivial splitting type? Is every uniform, admissible
locally-free sheaf with vanishing first Chern class trivial?

We've also seen that if $E$ is an admissible torsion-free and $\wp$ is a plane
in $\p3$, then the restriction $E|_{\wp}$ satisfies the following cohomological
condition:
$$ H^0(\wp,E|_{\wp}(k)) = 0 ~ \forall k\leq-1 ~~ , ~~
H^2(\wp,E|_{\wp}(k)) = 0 ~ \forall k\geq-2 ~~ .$$
A sheaf on $\p2$ satisfying the above conditions are called {\em instanton sheaves},
and are very interesting on their own right, also being closely related to instantons.
The analysis of how the instanton condition compares with semistability and trivial
splitting type is work in progress \cite{J-p2}, but many of the results proved here
have their analogs for instanton sheaves in $\p2$.

In particular, it is shown in \cite{J-p2} that every instanton sheaf is the
cohomology of a special monad, and that every rank 2 torsion-free instanton
sheaf is semistable. We can then conjecture that if a (rank 2) torsion-free
sheaf $E$ on $\mathbb{P}^k$ is the cohomology of a special monad, then $E$ 
is semistable; this is true for $k=2,3$.

 

 \end{document}